\documentclass[12pt]{amsart}
\newtheorem{thm}{Theorem}[section]
\newtheorem{theorem}[thm]{Theorem}
\newtheorem{corollary}[thm]{Corollary}
\newtheorem{lemma}[thm]{Lemma}
\newtheorem{proposition}[thm]{Proposition}
\theoremstyle{definition}

\newtheorem{example}[thm]{Example}
\theoremstyle{remark}
\newtheorem{remark}[thm]{Remark}
\numberwithin{equation}{section}

\newcommand{\ga}{\gamma}
\newcommand{\de}{\delta}

\newcommand{\ve}{\varepsilon}
\newcommand{\z}{\zeta}

\newcommand{\x}{\xi}
\renewcommand{\th}{\theta}

\newcommand{\ka}{\kappa}
\newcommand{\la}{\lambda}
\newcommand{\si}{\sigma}

\newcommand{\f}{\phi}

\newcommand{\h}{\chi}

\newcommand{\De}{\Delta}

\newcommand{\Om}{\Omega}

\newcommand{\R}{\mathbb R}

\newcommand{\dpa}{\partial}
\newcommand{\dr}{\mathrm{d}}

\newcommand{\sign}{\mathrm{sign}\,}

\newcommand{\supp}{\mathrm{supp}\,}

\newcommand{\meas}{\mathrm{meas}\,}

\newcommand{\const}{\mathrm{const\,}}

\begin{document}

\title[Fourier Tauberian Theorems]
{Fourier Tauberian Theorems and Applications}
\author{Yu. Safarov}
\address{
Department of Mathematics, King's College,
Strand, London WC2R 2LS, UK}
\email{ysafarov\@mth.kcl.ac.uk}

\subjclass{26A48, 35P15}
\keywords{Fourier Tauberian Theorems, Laplacian,
counting and spectral functions}
\date{February 2000}

\maketitle

Let $F$ be a non-decreasing function
and $\rho$ is an appropriate test
function on the real line $\R$. Then, under
certain conditions on the Fourier transform of the
convolution $\,\rho*F\,$, one can estimate the difference
$\,F-\rho*F\,$.
Results of this type are called Fourier Tauberian theorems.

The Fourier Tauberian theorems have been used by many authors
for the study of spectral asymptotics of
elliptic differential operators, with $F$ being either the
counting function or the spectral function (see, for example,
\cite{L}, \cite{H1}, \cite{H2}, \cite{DG}, \cite{I1},
\cite{I2}, \cite{S}, \cite{SV}). The required
estimate for $\,F-\rho*F\,$ was obtained under the assumption
that the derivative $\,\rho*F'\,$ admits a sufficiently good
estimate.

In applications $F$ often depends on additional
parameters and we are interested in estimates
which are uniform with respect to these parameters.
Then one has to assume that the
estimate for $\,\rho*F'\,$ holds uniformly and to
take this into account when estimating $\,F-\rho*F\,$.
As a result, there have been produced a
number of Fourier Tauberian theorems designed for
the study of various parameter dependent
problems. This has been done, in particular, for
semi-classical asymptotics (see, for example,
\cite{PP}). Note that all the authors used
the same idea of proof which goes back to the papers
\cite{L} and \cite{H1}.

The main aim of this paper is to present a general
version of Fourier Tauberian theorem which does not
require any a priori estimates of $\rho*F'$.
Our estimates contain only convolutions
of $F$ and test functions (see Section 1). This enables
one to obtain results which are uniform with respect to
any parameters involved.

Our proof is very different from the usual one.
It leads to more general results and, at the same time,
allows one to evaluate constants appearing
in the estimates (Section 2). Therefore our Tauberian
theorems can be used not only for the study of asymptotics
but also for obtaining explicit estimates of the spectral
and counting functions.

In particular, in Section 3, applying our Tauberian
theorems and Berezin's inequality, we
prove a refined version of the Li--Yau estimate
for the counting function of the Dirichlet Laplacian
in an arbitrary domain of finite volume. Our inequality
implies the Li--Yau estimate itself and, along with
that, the results on the asymptotic behaviour of the
counting function which are obtained by the
variational method. Note that the proof of the
Berezin inequality does not use variational
techniques. This implies that, even in the non-smooth case,
the classical asymptotic formulae can be proved without
referring to the Whitney decompositions and Dirichlet--Neumann
bracketing.

Throughout the paper
$\h_+$, $\h_-$ denote the characteristic functions
of the positive and negative semi-axes,
$\,\hat f(t)
:=(2\pi)^{-1/2}\int e^{-it\tau} f(\tau)\,\dr\tau\,$
is the Fourier transform of $f$, and
$\,\langle\tau\rangle:=\sqrt{1+\tau^2}\,$.

\section{Tauberian theorems I: basic estimates}
Let $F$ be a non-decreasing function on $\R$.
For the sake of definiteness,
we shall always be assuming that
\begin{equation}\label{1.1}
F(\tau)\ =\ \frac12\,\left[F(\tau+0)+F(\tau-0)\right]\,,
\qquad\forall\tau\in\R\,.
\end{equation}

\subsection{Auxiliary functions}
We shall deal with continuous functions
$\,\rho\,$ on $\R$ satisfying the following conditions:
\begin{enumerate}
\item[(1$_m$)]
$\,|\rho(\tau)|\le\const\langle\tau\rangle^{-2m-2}\,$,
where $m>-\frac12$;
\item[(2)\phantom{$_m$}]
$\,c_{\rho,0}:=\int\rho(\tau)\,\dr\tau=1\,$;
\item[(3)\phantom{$_m$}]
$\rho$ is even;
\item[(4)\phantom{$_m$}]
$\,\rho\ge0\,$;
\item[(5)\phantom{$_m$}]
$\,\supp\hat\rho\subset[-1,1]\,$.
\end{enumerate}
For every $m$ the functions $\rho$ satisfying (1$_m$)--(5)
do exist (see, for example, \cite{H2}, Section 17.5,
or Example \ref{E1.1} below).

\begin{example}\label{E1.1}
Let $l$ be a positive integer and
\begin{equation}\label{1.2}
\begin{array}{c}
\ga(\tau)\ :=\
\int_{-\frac\pi2}^{\frac\pi2}\,
(\frac\tau{2l}+s)^{-2l}\,
\sin^{2l}(\frac\tau{2l}+s)\,\dr s\,.
\end{array}
\end{equation}
The function $\ga$ satisfies (3), (5), and
\begin{equation}\label{1.3}
c_\ga^-\langle\tau\rangle^{-2l}\ \le\ \ga(\tau)
\ \le\ c_\ga^+\langle\tau\rangle^{-2l}
\end{equation}
with some positive constants $c_\ga^\pm$.
Indeed, (3) and (\ref{1.3}) are
obvious, and (5) follows from
the fact that $-2(2\pi)^{-1/2}\tau^{-1}\sin\tau$
is the Fourier transform of the characteristic
function of the interval $[-1,1]$. If
$\,\rho(\tau):=c_{\ga,0}^{-1}\ga(\tau)\,$ then
the conditions (2)--(5) are fulfilled and
(1$_m$) holds with $m=l-1$.
\end{example}

We shall always be assuming (1$_m$). Let
$$
\rho_{1,1}(\tau)\ :=\
\begin{cases}
\int_\tau^\infty\rho(\mu)\,\dr\mu,&\tau>0,\\
\phantom{-}0,&\tau=0,\\
-\int_{-\infty}^\tau\rho(\mu)\,\dr\mu,&\tau<0,
\end{cases}
$$
and, if (1$_m$) holds with $m>0$,
$$
\rho_{1,0}(\tau)\ :=\ \int_\tau^\infty\mu\,
\rho(\mu)\,\dr\mu\,,\qquad
\rho_{1,2}(\tau)\ :=\
\begin{cases}
\int_\tau^\infty\int_\mu^\infty
\rho(\la)\,\dr\la\,\dr\mu,&\tau\ge0,\\
\int_{-\infty}^\tau\int_{-\infty}^\mu
\rho(\la)\,\dr\la\,\dr\mu,&\tau\le0.
\end{cases}
$$
One can easily see that
$$
\rho_{1,0}(\tau)\le\const\langle\tau\rangle^{-2m},\quad
\rho_{1,1}(\tau)\le\const\langle\tau\rangle^{-2m-1},\quad
\rho_{1,2}(\tau)\le\const\langle\tau\rangle^{-2m}
$$
for all $\tau\ge0$.
Integrating by parts, we obtain
\begin{equation}\label{1.4}
\rho_{1,0}(\tau)\ =\ -\int_\tau^\infty\mu\,
\rho'_{1,1}(\mu)\,\dr\mu
\ =\ \rho_{1,2}(\tau)+\tau\,\rho_{1,1}(\tau)\,,
\qquad\forall\tau\ge0\,.
\end{equation}

Denote
\begin{equation*}
c_{\rho,\ka}\ :=\ \int|\mu|^\ka\,\rho(\mu)\,\dr\mu\,,
\qquad\forall\ka\in(-1,2m+1)\,.
\end{equation*}
Under condition (2), by Jensen's inequality, we have
\begin{equation}\label{1.5}
c_{\rho,r}^\ka\ \le\ c_{\rho,\ka}^r\,,
\qquad\forall\ka\ge r\ge0\,.
\end{equation}

If the condition (3) is fulfilled then
$\,\rho_{1,0}\,$ and $\,\rho_{1,2}\,$ are
even continuous functions, $\,\rho_{1,1}\,$
is an odd function continuous outside the origin
and
\begin{equation}\label{1.6}
\begin{array}{c}
\rho_{1,1}(\pm0)=\pm\,\frac12\,c_{\rho,0}\,,\qquad
\rho_{1,0}(0)=\rho_{1,2}(0)=\,\frac12\,c_{\rho,1}\,.
\end{array}
\end{equation}
Indeed, the first two equalities in (\ref{1.6})
are obvious, and the last follows from (\ref{1.4}).

The condition (4) and (\ref{1.4}) imply that
\begin{equation}\label{1.7}
\begin{array}{lr}
0\ \le\ \rho_{1,2}(\tau)\ \le\
\rho_{1,0}(\tau)\,,
&\qquad\forall\tau\ge0\,,\\
0\ \le\ \rho_{1,k}(\mu)\ \le\ \rho_{1,k}(\tau)\,,
&\qquad k=0,1,2,\quad\forall\mu\ge\tau\ge0\,.
\end{array}
\end{equation}

Let
\begin{equation}\label{1.8}
\rho_\de(\tau)\ :=\ \de\rho(\de\tau)\,,\qquad
\rho_{\de,k}(\tau)\ :=\ \de^{1-k}\,\rho_{1,k}(\de\tau)\,,
\qquad k=0,1,2,
\end{equation}
where $\de$ is an arbitrary positive number.
If (5) is fulfilled then
\begin{equation}\label{1.9}
\supp\hat\rho_{\de,0}\;\subset\;\supp\hat\rho_\de
\;\subset\;[-\de,\de]\,.
\end{equation}
Indeed, these inclusions follow from (\ref{1.8})
and the fact that $\,\rho_{1,0}\,$
is the convolution of the functions
$\,\mu\,\rho(\mu)\,$ and $\,\h_-(\mu)$.

\subsection{Main estimates}
If $f$ is a piecewise
continuous function on $\R^1$, we denote
\begin{eqnarray*}
f*F(\tau)&:=&\lim_{R\to\infty}
\int_{-R}^Rf(\tau-\mu)\,F(\mu)\,\dr\mu\,,\\
f*F'(\tau)&:=&\lim_{R\to\infty}
\int_{(-R,R)}f(\tau-\mu)\,\dr F(\mu)\,,
\end{eqnarray*}
whenever the limits exist.
We shall deduce the estimates for $\,F(\tau)\,$
from the following simple lemma.

\begin{lemma}\label{L1.2}
Let $\rho$ satisfy the conditions {\rm (1$_m$)--(3)} and
$\,\rho_{T,1}(\tau-s)\,F(s)\to0\,$ as
$\,s\to\pm\infty\,$ for some $\,T>0\,$ and 
$\,\tau\in\R\,$.  Then $\rho_{T,1}*F'(\tau)$
is well defined if and only if
$\,\rho_T*F(\tau)\,$ is well defined, and
\begin{equation}\label{1.10}
F(\tau)-\rho_T*F(\tau)\ =\ \rho_{T,1}*F'(\tau)\,.
\end{equation}
\end{lemma}

\begin{proof}
Integrating by parts, we obtain
\begin{multline*}
\int_{(-R,R)}\rho_{T,1}(\tau-\mu)\,\dr F(\mu)
=\int_{(-R,\tau)}\rho_{T,1}(\tau-\mu)\,\dr F(\mu)+
\int_{(\tau,R)}\rho_{T,1}(\tau-\mu)\,\dr F(\mu)\\
=\ -\int_{-R}^R\rho_T(\tau-\mu)\,F(\mu)\,\dr\mu
+\rho_{T,1}(+0)\,F(\tau-0)-\rho_{T,1}(-0)\,F(\tau+0)\\
-\rho_{T,1}(\tau+R)\,F(-R+0)+\rho_{T,1}(\tau-R)\,F(R-0)\,.
\end{multline*}
In view of (\ref{1.1}), (\ref{1.6}) and (2),
we have
$$
\rho_{T,1}(+0)\,F(\tau-0)-
\rho_{T,1}(-0)\,F(\tau+0)=F(\tau)\,.
$$
Now the lemma is proved by passing to the limit
as $R\to\infty$.
\end{proof}

\begin{theorem}\label{T1.3}
Let $\rho$ satisfy the conditions 
{\rm (1$_m$)--(4)} with $m>0\,$. Assume that 
$\,\rho_{\de,0}(\tau-s)\,F(s)\to0\,$ as
$\,s\to\pm\infty\,$ 
and $\,\rho_{\de,0}*F'(\tau)<\infty\,$
for some $\,\de>0\,$ and $\,\tau\in\R\,$.
Then $\,\rho_T*F(\tau)<\infty\,$ and
\begin{equation}\label{1.11}
|F(\tau)-\rho_T*F(\tau)|\ \le\
c_{\rho,1}^{-1}\,\de^{-1}\,\rho_{\de,0}*F'(\tau)
\end{equation}
for all $\,T\ge\de$.
\end{theorem}

\begin{proof}
The identity (\ref{1.4}) and (4) imply that
$$
\frac{\dr\hfill}{\dr\tau}\left(
\frac{\rho_{1,1}(\tau)}{\rho_{1,0}(\tau)}\right)
\ =\ \frac{\rho(\tau)
\left(\tau\,\rho_{1,1}(\tau)-\rho_{1,0}(\tau)\right)}
{\left(\rho_{1,0}(\tau)\right)^2}\ \le\ 0\,,
\qquad\forall\tau>0\,.
$$
Therefore, in view of (2) and (\ref{1.6}),
\begin{equation*}
\frac{|\rho_{1,1}(\tau)|}{\rho_{1,0}(\tau)}\ \le\
\frac{|\rho_{1,1}(+0)|}{\rho_{1,0}(0)}
\ =\ c_{\rho,1}^{-1}\,,\qquad\forall\tau>0\,.
\end{equation*}
Taking into account (3), (\ref{1.8}) and the 
second inequality (\ref{1.7}), we obtain
\begin{equation}\label{1.12}
|\rho_{T,1}(\tau)|\ \le\
\frac{\rho_{T,0}(\tau)}{c_{\rho,1}\,T}\ \le\
\frac{\rho_{\de,0}(\tau)}{c_{\rho,1}\,\de}\,,
\qquad\forall T\ge\de>0,
\quad\forall\tau\in\R\,.
\end{equation}

The inequality 
(\ref{1.12}) implies that $\rho$ and $F$
satisfy the conditions of Lemma \ref{L1.2} and that
$\,\rho_T*F(\tau)<\infty\,$. Obviously, (\ref{1.11})
follows from (\ref{1.10}) and (\ref{1.12}).
\end{proof}

\begin{remark}\label{R1.4}
If $T=\de$ then the estimate (\ref{1.11}) can be rewritten in the
form
\begin{equation}\label{1.13}
\rho_\de^+*F(\tau)\ \le\ F(\tau)\ \le\ \rho_\de^-*F(\tau)\,,
\end{equation}
where
$\,
\rho_\de^{\pm}(\tau)
:=\rho_\de(\tau)\pm c_{\rho,1}^{-1}\,\de\tau\,\rho_\de(\tau)
\,$.
\end{remark}

\begin{remark}\label{R1.5}
The inequality (\ref{1.11})
is not precise in the sense that, apart from some
degenerate situations, it never turns into an equality.
The crucial point in our proof is the estimate
$\,|\rho_{T,1}|\le c_{\rho,1}^{-1}\,\de^{-1}\,\rho_{\de,0}\,$
which implies that
$\,|\rho_{T,1}*F'|\le
c_{\rho,1}^{-1}\,\de^{-1}\,\rho_{\de,0}*F'(\tau)\,$. However,
the function $\rho_{T,1}$ is negative on one half-line
and positive on another, so $\,|\rho_{T,1}*F'|\,$
may well admit much a better estimate.
Using this observation, one can try to improve our
results under additional conditions on the function
$F$.
\end{remark}

\begin{theorem}\label{T1.6}
Let $\,[a,b]\,$ be a bounded interval.
Assume that the conditions of Theorem \ref{T1.3}
are fulfilled for every $\,\tau\in[a,b]\,$
and that
$\,\rho_{\de,0}*F'(\tau)\,$ is uniformly bounded 
on $\,[a,b]\,$. Then 
\begin{eqnarray}
&&-\,T^{-1}\de^{-1}\,f(b)\;\rho_{\de,0}*F'(b)\nonumber\\ 
&&\le\ \int_a^b f(\tau)
\left[F(\tau)-\rho_T*F(\tau)\right]\dr\tau\label{1.14}\\
&&\le\ T^{-1}\de^{-1}\,f(a)\;\rho_{\de,0}*F'(a)
+T^{-1}\de^{-1}\int_a^b f'(\tau)\,\rho_{\de,0}*F'(\tau)\,\dr\tau\nonumber
\end{eqnarray}
for every non-negative non-decreasing function
$\,f\in C^1[a,b]\,$ and all $T\ge\de$.
\end{theorem}

\begin{proof}
In view of (\ref{1.7}) and (\ref{1.8}) we have
\begin{equation}\label{1.15}
T\rho_{T,2}(\tau)\ \le\
T^{-1}\rho_{T,0}(\tau)\ \le\
\de^{-1}\rho_{\de,0}(\tau)\,,
\qquad\forall T\ge\de>0\,,
\quad\forall\tau\in\R\,.
\end{equation}
This estimates, (\ref{1.12}) and Lemma \ref{L1.2}
imply that the functions $\,\rho_{T,2}*F'(\tau)\,$,
$\,|\rho_{T,1}|*F'(\tau)\,$ and
$\,\rho_T*F(\tau)\,$
are uniformly bounded on $\,[a,b]\,$.
Since $\,\rho'_{T,2}(s)=-\rho_{T,1}(s)\,$
whenever $\,s\ne0\,$, integrating by parts
with respect to $\tau$ we obtain
\begin{multline*}
\int_a^b f(\tau)\int\rho_{T,1}(\tau-\mu)\,
\dr F(\mu)\,\dr\tau\ =\
f(a)\int\rho_{T,2}(a-\mu)\,\dr F(\mu)\\
-f(b)\int\rho_{T,2}(b-\mu)\,\dr F(\mu)
+\int_a^b f'(\tau)\left(\int\rho_{T,2}(\tau-\mu)\,
\dr F(\mu)\right)\dr\tau\,.
\end{multline*}
Now (\ref{1.14}) follows from  Lemma \ref{L1.2}
and (\ref{1.15}).
\end{proof}

If $f\equiv1$ then (\ref{1.14}) turns into
\begin{equation}\label{1.16}
-\,\rho_{\de,0}*F'(b)\ \le\
T\de\int_a^b
\left[F(\mu)-\rho_T*F(\mu)\right]\dr\mu
\ \le\ \rho_{\de,0}*F'(a)\,.
\end{equation}
This estimate and
the obvious inequalities
\begin{equation}\label{1.17}
\ve^{-1}\int_{\tau-\ve}^\tau F(\mu)\,\dr\mu
\ \le\ F(\tau)\ \le\
\ve^{-1}\int_\tau^{\tau+\ve}F(\mu)\,\dr\mu\,,
\qquad\forall\ve>0\,,
\end{equation}
imply the following corollary.

\begin{corollary}\label{C1.7}
Under conditions of Theorem \ref{T1.6}
\begin{eqnarray}
F(b)&\ge&
\ve^{-1}\int_{b-\ve}^b\rho_T*F(\mu)\,\dr\mu
-\,\ve^{-1}T^{-1}\de^{-1}\,\rho_{\de,0}*F'(b)\,,
\label{1.18}\\
F(a)&\le&
\ve^{-1}\int_a^{a+\ve}\rho_T*F(\mu)\,\dr\mu
+\,\ve^{-1}T^{-1}\de^{-1}\,\rho_{\de,0}*F'(a)
\label{1.19}\end{eqnarray}
for all $\,\ve\in(0,b-a]\,$ and $\,T\ge\de\,$.
\end{corollary}

If (4) is fulfilled then $\,\rho_T*F\,$ is a non-decreasing
function. Therefore (\ref{1.18}) and (\ref{1.19}) imply that
\begin{eqnarray}
F(b)&\ge&\rho_T*F(b-\ve)
-\ve^{-1}T^{-1}\de^{-1}\,\rho_{\de,0}*F'(b)\,,
\label{1.20}\\
F(a)&\le&\rho_T*F(a+\ve)
+\ve^{-1}T^{-1}\de^{-1}\,\rho_{\de,0}*F'(a)\,.
\label{1.21}\end{eqnarray}

\begin{remark}\label{R1.8}
It is clear from the proof that Theorems
\ref{T1.3} and \ref{T1.6} remain
valid (with some other constants independent 
of $\de$ and $T$)
if we drop the condition (4) and
replace $\rho_{\de,0}(\tau)$ with
an arbitrary non-negative function
$\ga_\de$ such that
$\,|\rho_{T,1}(\tau)|\le\const\,
\de^{-1}\ga_\de(\tau)\,$ and
$\,|\rho_{T,2}(\tau)|\le\const\,
T^{-1}\de^{-1}\ga_\de(\tau)\,$.
In particular, one can take 
$\ga_{\de}(\tau)=\de\ga(\de\tau)$, where $\ga$
is the function defined by (\ref{1.3}) with
$\,l=m\,$.
\end{remark}

\section{Tauberian theorems II: applications}

\subsection{General remarks}
From now on we shall be assuming that the function
$F$ is polynomially bounded. Then
the conditions of Theorems \ref{T1.3} and \ref{T1.6}
are fulfilled for all $\tau,a,b\in\R^1$ and
$T\ge\de>0$ whenever $\rho$ satisfies (1$_m$)
with a sufficiently large $m$.

So far we have not used the condition (5),
which is not needed to prove the estimates.
However, this condition often appears
in applications. It implies that the convolutions
$\rho_T*F$ and $\rho_{T,0}*F'$
are determined by the restrictions of $\hat F$
to the interval $(-T,T)$. If
\begin{equation}\label{2.1}
\left.\hat F_0(t)\right|_{(-T,T)}
=\ \left.\hat F(t)\right|_{(-T,T)}
\end{equation}
then, under condition (5),
$\rho_T*F=\rho_T*F_0$
and $\rho_{\de,0}*F'=\rho_{\de,0}*F'_0$
for all $\de\le T$. If $F_0(\tau)$
behaves like a linear combination of homogeneous
functions for large $\tau$ then $\rho_{\de,0}*F'_0$
is of lower order than $\rho_T*F_0$, so it plays the
role of an error term in asymptotic formulae.

It is not always possible to find a model function
$F_0$ satisfying (\ref{2.1}). However, one can
often construct $\tilde F_0$
in such a way that the convolutions
$\rho_T*(F-\tilde F_0)(\tau)$
and $\rho_{\de,0}*(F'-\tilde F'_0)(\tau)$
admit good estimates for large $\tau$
(roughly speaking, it happens if the Fourier
transforms of $F$ and $\tilde F_0$ have similar
singularities on the corresponding interval).
Then the Tauberian theorems imply estimates
with the error term
$$
\pm\left(|\rho_T*(F-\tilde F_0)(\tau)|
+|\rho_{\de,0}*(F'-\tilde F'_0)(\tau)|\right).
$$
In particular, if $F$ is the spectral or counting
function of an elliptic partial differential operator
with smooth coefficients
then (\ref{1.11}) gives a precise reminder
estimate in the Weyl asymptotic formula, and
the refined estimates (\ref{1.20}), (\ref{1.21}) allow
one to obtain the second asymptotic term
by letting $T\to\infty$ (see \cite{SV} for details).

In applications to the second order differential
operators it is usually more convenient to deal with
the cosine Fourier transform of $F'$.
The following elementary observation enables one
to apply our results in the case where information
on the sine Fourier transform of $F'$ is not available.

\begin{proposition}\label{P2.1}
If the cosine Fourier transforms of the derivatives
$F'$ and $F'_0$ coincide on an
interval $(-\de,\de)$ then the Fourier transforms of
the functions $F(\tau)-F(-\tau)$ and $F_0(\tau)-F_0(-\tau)$
coincide on the same interval.
\end{proposition}

\subsection{Test functions $\rho$}
In this subsection we consider a class of
functions $\rho$ satisfying (1$_m$)--(5)
and estimate the constants $c_{\rho,\ka}$.

\begin{lemma}\label{L2.2}
Let $\,\z\in C^{m+1}[-\frac12,\frac12]\,$
be a real-valued even function such that
$\,\|\z\|_{L_2}=1\,$ and
$\,\z^{(k)}(\pm\frac12)=0\,$
for $k=0,1,\ldots m-1\,$, where
$\z^{(k)}$ denotes the $k$th derivative.
If we extend $\z$ to $\R$ by zero then
$\,\rho:=(\hat\z)^2\,$ satisfies
{\rm (1$_m$)--(5)} and
\begin{equation}\label{2.2}
c_{\rho,2k}\ =\ \|\z^{(k)}\|_{L_2}^2\,,
\qquad k=0,1,\ldots,m\,.
\end{equation}
\end{lemma}

\begin{proof}
The conditions (3) and (4) are
obviously fulfilled; (2),
(5) and (\ref{2.2}) follow from the fact that
$\,\hat\rho=(2\pi)^{-1/2}\,\z*\z\,$. Finally,
(1$_m$) holds true
because the $(m+1)$th derivative of the extended
function $\z$ coincides with a linear combination of
an $L_1$-function and two $\de$-functions.
\end{proof}

The following lemma is a consequence of
the uncertainty principle. 

\begin{lemma}\label{L2.3}
If $\rho$ is defined as in Lemma \ref{L2.2} then
\begin{equation}\label{2.3}
\begin{array}{c}
c_{\rho,1}\ \ge\ \frac\pi2\,.
\end{array}
\end{equation}
\end{lemma}

\begin{proof}
Let $\Pi_a$ be the
multiplication operator and $\hat\Pi_a$ be
the Fourier multiplier generated by
the characteristic function of the interval $[-a,a]$.
Then the Hilbert-Schmidt norm of the operator
$\,\hat\Pi_{a_1}\Pi_{a_2}\,$ acting in $L_2(\R)$ is equal to
$\,\sqrt{2\pi^{-1}a_1a_2}\,$. Therefore
\begin{equation*}
2\int_0^\mu\hat\z^2(\tau)\,\dr\tau
\ =\ \|\hat\Pi_\mu\Pi_{1/2}\z\|_{L_2}^2
\ \le\ \pi^{-1}\mu\,\|\z\|_{L_2}^2
\ =\ \pi^{-1}\mu\,,
\end{equation*}
which implies that
\begin{multline*}
c_{\rho,1}\ =\
2\int_0^\infty\mu\,\hat\z^2(\mu)\,\dr\mu
\ =\ 2\int_0^\infty\int_\mu^\infty
\hat\z^2(\tau)\,\dr\tau\,\dr\mu\\
\ge\ 2\int_0^\pi\int_\mu^\infty
\hat\z^2(\tau)\,\dr\tau\,\dr\mu
\ \ge\ \int_0^\pi(1-\pi^{-1}\mu)\,\dr\mu\ =\ \frac\pi2\,.
\end{multline*}
\end{proof}

\begin{remark}\label{R2.4}
As follows from Nazarov's theorem (see \cite{Na}
or \cite{HJ}),
$$
\begin{array}{c}
\int_\mu^\infty\hat\f^2(\tau)\,\dr\tau
\ \ge\ b_1e^{-b_2\mu}\,,
\qquad\forall\f\in C_0^\infty(-\frac12,\frac12)\,,
\quad\forall\mu\ge0\,,
\end{array}
$$
where $b_1,b_2>0$ are some absolute constants.
Using the estimates for $b_1,b_2$ obtained in 
\cite{Na}, one can slightly improve
the estimate (\ref{2.3}).
\end{remark}

\begin{example}\label{E2.5}
Let $\tilde\nu_m$ be the first eigenvalue of
the operator $\displaystyle\frac{\dr^{2m}\hfill}{\dr t^{2m}}$
on the interval $(-\frac12,\frac12)$ subject to Dirichlet
boundary condition, and let $\z_m$ be the corresponding
real even normalized eigenfunction. Denote
$\displaystyle\nu_m:=\left(\tilde\nu_m\right)^{\frac1{2m}}$.
If we define $\rho$ as in Lemma \ref{L2.2}
then, in view of (\ref{2.2}) and (\ref{1.5}),
\begin{equation}\label{2.4}
c_{\rho,2m}\ =\ \nu_m^{2m}\,,\qquad\qquad
c_{\rho,\ka}\ \le\ \nu_m^\ka\,,\qquad\forall\ka<2m\,.
\end{equation}
\end{example}

The eigenvalues $\tilde\nu_m=\nu_m^{2m}$ grow very fast as
$m\to\infty$. The following lemma
gives a rough estimate for $\nu_m$.

\begin{lemma}\label{L2.6}
We have $\;\nu_m\le2m\sqrt[2m]3\;$ for all $m\ge2$.
\end{lemma}

\begin{proof}
If $\,\f(t)=\left(\frac14-t^2\right)^m\,$
and $\|\cdot\|_{L_2}$ is the norm
in $L_2\left(-\frac12,\frac12\right)$
then
\begin{equation}\label{2.5}
\tilde\nu_m\ \le\
\frac{\|\f^{(m)}\|_{L_2}^2}{\|\f\|_{L_2}^2}
\ =\ \frac{(4m+1)!\,(m!)^2}{(2m+1)!\,(2m)!}
\ \le\ 2^{2m+1}\,(2m)!\,.
\end{equation}
One can easily see that
$$
 \frac{2^{2m}\,(2m)!}{(2m)^{2m}}\ =\
\frac{2\,(m^2-1)\dots(m^2-(m-1)^2)}{m^{2m-2}}
\ \le\ \frac{2\,(m^2-(m-1)^2)}{m^2}
\ \le\ \frac32\,.
$$
Therefore (\ref{2.5}) implies the required estimate.
\end{proof}

\subsection{Power like singularities}
Assume that $|F(\tau)|\le\const\,(|\tau|+1)^n$
with a non-negative integer $n$ and define
$$
\si_n:=\begin{cases}
0\,,&\text{if $n$ is odd,}\\
1\,,&\text{if $n$ is even,}
\end{cases}
\qquad m_n:=\begin{cases}
\frac{n+1}2\,,&\text{if $n$ is odd,}\\
\frac{n+2}2\,,&\text{if $n$ is even,}
\end{cases}
$$
$$
P_n^+(\tau,\mu):=\frac{(\tau+\mu)^n
+(\tau-\mu)^n}2\,,\qquad
P_n^-(\tau,\mu):=\frac{\mu\,(\tau+\mu)^n
-\mu\,(\tau-\mu)^n}2\;.
$$
Clearly, $P_n^\pm$ are homogeneous polynomials
in $(\tau,\mu)$ with positive coefficients, which
contain only even powers of $\mu$.

\begin{lemma}\label{L2.7}
Let $\rho$ be a function satisfying 
{\rm (3)}, {\rm (5)} and {\rm (1$_m$)} with $m>\frac{n}2\,$.
If $\,\supp F\subset(0,+\infty)\,$ and the cosine Fourier 
transform of $\,F'(\tau)\,$ coincides on the interval 
$\,(-\de,\de)\,$ with the cosine Fourier transform of 
the function $\,n\tau_+^{n-1}\,$ then
\begin{eqnarray}
\rho_\de*F(\tau)&\ge&
\int\left[P_n^+(\tau,\de^{-1}\mu)-
\si_n\,\de^{-n}|\mu|^n\right] \rho(\mu)\,\dr\mu\,,
\label{2.6}\\
\rho_\de*F(\tau)&\le&
\int P_n^+(\tau,\de^{-1}\mu)\,\rho(\mu)\,\dr\mu\,,
\label{2.7}\\
\rho_{\de,0}*F'(\tau)&\le&
\de^2\int\left[P_n^-(\tau,\de^{-1}\mu)
+\si_n\,\de^{-n-1}|\mu|^{n+1}\right]\rho(\mu)\,\dr\mu
\label{2.8}
\end{eqnarray}
for all $\tau>0$.
\end{lemma}

\begin{proof}
According to Proposition \ref{P2.1}, the Fourier
transform of $F(\tau)-F(-\tau)$ coincides on the
interval $(-\de,\de)$ with the Fourier transform
of
$$
\sign\tau\,|\tau|^n\ =\ \left(1
-2\si_n\,\h_-(\tau)\right)\tau^n\,.
$$
Since $\rho$ is even, this implies that
\begin{multline*}
\rho_\de*F(\tau)\ =\ \de\int
\left(1-2\si_n\,\h_-(\tau-\mu)\right)(\tau-\mu)^n\,
\rho(\de\mu)\,\dr\mu\\
=\ \int P_n^+(\tau,\de^{-1}\mu)\,\rho(\mu)\,\dr\mu
\;-\;2\si_n\int_{\de\tau}^\infty
(\de^{-1}\mu-\tau)^n\,\rho(\mu)\,\dr\mu\,,
\end{multline*}
\begin{multline*}
\rho_{\de,0}*F'(\tau)=\rho'_{\de,0}*F(\tau)
=-\,\de^3\int\left(1-2\si_n\,\h_-(\tau-\mu)\right)
(\tau-\mu)^n\,\mu\,\rho(\de\mu)\,\dr\mu\\
=\ \de^2\int P_n^-(\tau,\de^{-1}\mu)\,\rho(\mu)\,\dr\mu
\;+\;2\si_n\,\de\int_{\de\tau}^\infty
(\de^{-1}\mu-\tau)^n\,\mu\,\rho(\mu)\,\dr\mu
\end{multline*}
for all $\tau>0$. Estimating
$\,0\le(\de^{-1}\mu-\tau)\le\de^{-1}\mu\,$
in the integrals on the right hand sides,
we arrive at (\ref{2.6})--(\ref{2.8}).
\end{proof}

The obvious inequalities
\begin{eqnarray*}
&\tau^n+\si_n\,|\nu|^n\ \le\ P_n^+(\tau,\nu)
\ \le\
\tau^n+n\,|\nu|\,(\tau+|\nu|)^{n-1}\,,\\
&P_n^-(\tau,\nu)+\si_n\,|\nu|^{n+1}\ \le\
n\,\nu^2\,(\tau+|\nu|)^{n-1}
\end{eqnarray*}
and (\ref{2.6})--(\ref{2.8}) imply that, for all
$\tau>0$,
\begin{eqnarray}
0\ \le\ \rho_\de*F(\tau)-\tau^n
&\le&
n\,\de^{-1}\int|\mu|\,(\tau+\de^{-1}|\mu|)^{n-1}\,
\rho(\mu)\,\dr\mu\,,\label{2.9}\\
\rho_{\de,0}*F'(\tau)
&\le&
n\int\mu^2\,(\tau+\de^{-1}|\mu|)^{n-1}\,\rho(\mu)\,\dr\mu\,.
\label{2.10}
\end{eqnarray}

Note that $m_n$ is the minimal positive integer
which is greater than $\frac{n}2$.
If $\rho$ is defined as in Lemma \ref{L2.2}
with $m=m_n$ then, by (\ref{2.2}),
\begin{equation}\label{2.11}
\int P_n^\pm(\tau,\de^{-1}\mu)\,\rho(\mu)\,\dr\mu
\ =\ \left(P_n^\pm(\tau,\de^{-1}D_t)\z,\z\right)_{L_2}.
\end{equation}
Applying (\ref{2.6})--(\ref{2.11}) and (\ref{1.11}) or
(\ref{1.18}), (\ref{1.19}),
one can obtain various estimates for $F(\tau)$.

\begin{example}\label{E2.8}
Let $n=3$ and $\z$ be an arbitrary function
satisfying conditions of Lemma \ref{L2.2} with
$m=m_n=2$. If the conditions of Lemma \ref{L2.7}
are fulfilled then (\ref{2.6})--(\ref{2.8}), (\ref{2.11})
and (\ref{1.19}), (\ref{1.20}) with $T=\de$
imply that
\begin{eqnarray*}
F(\tau)&\ge&
\tau^3-\frac{3\ve\tau^2}2+\ve^2\tau-\frac{\ve^3}4
+\frac3{2\de^2}\left(\tau-\frac{\tau^2}\ve-\frac\ve2\right)\|\z'\|_{L_2}^2
-\frac1{\ve\de^4}\,\|\z''\|_{L_2}^2\,,\\
F(\tau)&\le&
\tau^3+\frac{3\ve\tau^2}2+\ve^2\tau+\frac{\ve^3}4
+\frac3{2\de^2}\left(\tau+\frac{\tau^2}\ve+\frac\ve2\right)\|\z'\|_{L_2}^2
+\frac1{\ve\de^4}\,\|\z''\|_{L_2}^2
\end{eqnarray*}
for all $\ve>0$ and $\tau>0$. Thus,
$F(\tau)$ lies between the first Dirichlet eigenvalues
of ordinary differential operators
generated by the quadratic forms on the right hand sides of
the above inequalities.
\end{example}

\begin{corollary}\label{C2.9}
Under conditions of Lemma \ref{L2.7}
\begin{eqnarray}
F(\tau)&\ge&
\tau^n-2\pi^{-1}\nu_{m_n}^2\,n\,\de^{-1}
(\tau+\de^{-1}\nu_{m_n})^{n-1},\label{2.12}\\
F(\tau)&\le&
\tau^n+(2\pi^{-1}\nu_{m_n}^2+\nu_{m_n})\,n\,\de^{-1}\,
(\tau+\de^{-1}\nu_{m_n})^{n-1}\label{2.13}
\end{eqnarray}
for all $\tau>0$.
\end{corollary}

\begin{proof}
If we define $\rho$ as in Lemma \ref{L2.2}
with $\z=\z_m$ (see Example \ref{E2.5}) then
(\ref{2.12}), (\ref{2.13}) follow from (\ref{1.11})
with $T=\de$,
(\ref{2.9}), (\ref{2.10}), (\ref{2.3}) and (\ref{2.4}).
\end{proof}

\begin{corollary}\label{C2.10}
Under conditions of Lemma \ref{L2.7}
\begin{eqnarray}
\int_0^{\la^2} F(\sqrt\mu)\,\dr\mu
&\ge&
\frac{2\,\la^{n+2}}{n+2}\,
-2n\,\nu_{m_n}^2\de^{-2}\,\la\,
(\la+\de^{-1}\nu_{m_n})^{n-1}\,,\label{2.14}\\
\int_0^{\la^2} F(\sqrt\mu)\,\dr\mu
&\le&
\frac{2\,\la^{n+2}}{n+2}\,+\,
(n+1)\,\nu_{m_n}^2\de^{-2}\,(\la+\de^{-1}\nu_{m_n})^n
\label{2.15}
\end{eqnarray}
for all $\la>0$.
\end{corollary}

\begin{proof}
Since $\,\int_0^{\la^2}F(\sqrt\mu)\,\dr\mu=
2\int_0^\la F(\tau)\,\tau\,\dr\tau\,$, Theorem
\ref{T1.6} with $T=\de$, $a=0$, $b=\la$
and $f(\tau)=\tau$ implies
\begin{eqnarray}
\int_0^{\la^2}F(\sqrt\mu)\,\dr\mu &\ge&
2\int_0^\la\tau\,\rho_\de*F(\tau)\,\dr\tau
-2\de^{-2}\,\la\;\rho_{\de,0}*F'(\la)\,,\label{2.16}\\
\int_0^{\la^2}F(\sqrt\mu)\,\dr\mu &\le&
2\int_0^\la\left(\tau\,\rho_\de*F(\tau)+\de^{-2}
\rho_{\de,0}*F'(\tau)\right)\dr\tau\,.\label{2.17}
\end{eqnarray}

Let $\rho$ be defined as in Lemma \ref{L2.2}
with $\z=\z_m$. Then (\ref{2.14}) follows from
(\ref{2.16}), (\ref{2.9}), (\ref{2.10}) and (\ref{2.4}).
Since
$\,\tau\,P_n^+(\tau,\nu)+P_n^-(\tau,\nu)=P_{n+1}^+(\tau,\nu)\,$,
the inequality (\ref{2.17}) and (\ref{2.7}), (\ref{2.8})
imply that
$$
\int_0^{\la^2}F(\sqrt\mu)\,\dr\mu\ \le\
2\int_0^\la\int\left(P_{n+1}^+(\tau,\de^{-1}\mu)
+\si_n\,|\de^{-1}\mu|^{n+1}\right)\rho(\mu)\,\dr\mu\,
\dr\tau\,.
$$
Estimating
\begin{equation*}
\int_0^\la\left[P_{n+1}^+(\tau,\nu)
+\si_n\,|\nu|^{n+1}\right]\dr\tau
\ \le\ 
\frac{\la^{n+2}}{n+2}+\frac{n+1}2\,\nu^2\,(\la+|\nu|)^n
\end{equation*}
with $\,\nu=\de^{-1}\mu\,$
and applying (\ref{2.4}), we obtain (\ref{2.15}).
\end{proof}

\section{Applications to the Laplace operator}

Let $\Om\subset\R^n$ be an open domain and
$d(x)$ be the distance from $x\in\Om$
to the boundary $\dpa\Om$.

\subsection
{Estimates of the spectral function}
Consider the Laplacian $\De_B$ in $\Om$ subject
to a self-adjoint boundary condition
$\left.B(x,D_x)u\right|_{\dpa\Om}=0$,
where $B$ is a differential operator. Assume that
the operator $-\De_B$ is non-negative and denote by
$\Pi(\la)$ its spectral projection corresponding
to the interval $[0,\la)$. Let $e(x,y;\la)$
be the integral kernel of the operator
$\frac{\Pi(\la-0)+\Pi(\la+0)}2$
(the so-called {\sl spectral function}).
The Sobolev embedding theorem implies that
$e(x,y;\la)$ is a smooth function on $\Om\times\Om$
for each fixed $\la$ and that $e(x,x;\la)$ is a
non-decreasing polynomially bounded function of $\la$
for each fixed $x\in\Om$.

Let $\De_0$ be the Laplacian on $\R^n\,$,
and $\,e_0(x,y;\la)$, $\tilde e_0(x,y;\la)$,
$\tilde e(x,y;\la)$ be the spectral functions
of the operators $\De_0$, $\sqrt\De_0$, $\sqrt{\De_B}$
respectively. Then
\begin{eqnarray*}
\h_+(\tau)\,e(x,x;\tau^2)&=&\tilde e(x,x;\tau)\,,\\
\h_+(\tau)\,e_0(x,x;\tau^2)&=&\tilde e_0(x,x;\tau)
\ =\ C_n\,\tau_+^n\,,
\end{eqnarray*}
where
\begin{equation}
C_n\ :=\ (2\pi)^{-n}\,
\meas\{\x\in\R^n:|\x|<1\}\,.\label{3.1}
\end{equation}

By the spectral theorem, the cosine Fourier transform
of $\frac{\dr\hfill}{\dr\tau}\tilde e(x,y;\tau)$
coincides with the fundamental solution
$u(x,y;t)$ of the wave equation in $\Om$,
\begin{equation*}
u_{tt}=\De u\,,\qquad\left.Bu\right|_{\dpa\Om}=0\,,
\qquad\left.u\right|_{t=0}=\de(x-y)\,,
\qquad\left.u_t\right|_{t=0}=0\,.
\end{equation*}
Due to the finite speed of propagation, $u(x,x;t)$
is equal to $u_0(x,x;t)$ whenever
$t\in(-d(x),d(x))$, where $u_0(x,y;t)$ is
the fundamental solution of the wave equation in $\R^n$.
Thus, the cosine Fourier transforms
of the derivatives
$\frac{\dr\hfill}{\dr\tau}\tilde e_0(x,x;\tau)$
and $\frac{\dr\hfill}{\dr\tau}\tilde e(x,x;\tau)$
coincide on the time interval $(-d(x),d(x))$.
Applying (\ref{2.12})--(\ref{2.15}) to
$F(\tau)=C_n^{-1}\,\tilde e(x,x;\tau)$ we obtain
the following corollary.

\begin{corollary}\label{C3.1}
For every $x\in\Om$ and all $\la>0$ we have
\begin{eqnarray}
&&e(x,x;\la)\ \ge\
C_n\,\la^{n/2}\;-\;\frac{n\,C_n\,2\pi^{-1}\nu_{m_n}^2}{d(x)}\,
\left(\la^{1/2}+\frac{\nu_{m_n}}{d(x)}\right)^{n-1},
\label{3.2}\\
&&e(x,x;\la)\ \le\
C_n\,\la^{n/2}\;+\;
\frac{n\,C_n\,(2\pi^{-1}\nu_{m_n}^2+\nu_{m_n})}{d(x)}\,
\left(\la^{1/2}+\frac{\nu_{m_n}}{d(x)}\right)^{n-1},
\label{3.3}
\end{eqnarray}
\begin{eqnarray}
&&\int_0^\la e(x,x;\mu)\,\dr\mu\ \ge\
\frac{2\,C_n\,\la^{n/2+1}}{n+2}\,
\;-\;\frac{2n\,C_n\,\nu_{m_n}^2\,\la^{1/2}}{(d(x))^2}\,
\left(\la^{1/2}+\frac{\nu_{m_n}}{d(x)}\right)^{n-1}
\label{3.4}\\
&&\int_0^\la e(x,x;\mu)\,\dr\mu\ \le\
\frac{2\,C_n\,\la^{n/2+1}}{n+2}\,
\;+\;\frac{(n+1)\,C_n\,\nu_{m_n}^2}{(d(x))^2}\,
\left(\la^{1/2}+\frac{\nu_{m_n}}{d(x)}\right)^n.
\label{3.5}\end{eqnarray}
\end{corollary}

\subsection
{Estimates of the counting function of the Dirichlet
Laplacian}
In this subsection we shall be assuming that
$|\Om|<\infty$, where $|\cdot|$ denotes the
$n$-dimensional Lebesgue measure.

Consider the positive operator $-\De_D$, where $\De_D$
is the Dirichlet Laplacian in $\Om$. Let
$N(\la)$ be the number of its eigenvalues lying
below $\la$.
The following theorem is due to F. Berezin
\cite{B}.

\begin{theorem}\label{T3.2} For all $\la\ge0$
we have
\begin{equation}\label{3.6}
\int_0^\la N(\mu)\,\dr\mu
\ \le\
\frac{2}{n+2}\;C_n\,|\Om|\,\la^{n/2+1}\,.
\end{equation}
\end{theorem}

This results was reproduced in \cite{La}. A. Laptev
also noticed that the famous Li--Yau estimate
\begin{equation}\label{3.7}
N(\la)\ \le\ (1+2/n)^{n/2}\,C_n\,|\Om|\,\la^{n/2}\,,
\qquad\forall\la\ge0\,,
\end{equation}
(see \cite{LY})
is a one line consequence of (\ref{3.6}). Indeed,
(\ref{3.7}) can be proved by estimating
\begin{equation}\label{3.8}
N(\la)\ \le\
(\th\la)^{-1}\int_0^{\la+\th\la}N(\mu)\,\dr\mu
\ \le\
\frac{2\,(1+\th)^{n/2+1}}{(n+2)\,\th}\;
C_n\,|\Om|\,\la^{n/2}
\end{equation}
and optimizing the choice of $\th>0\,$.

\begin{remark}\label{R3.3}
In \cite{B} F. Berezin proved an analogue of (\ref{3.6})
for general operators with constant
coefficients subject to Dirichlet boundary condition.
In the same way as above, applying the first inequality
(\ref{3.8}) and Berezin's estimates, one can easily obtain
upper  bounds for the corresponding counting functions
(see \cite{La}).
\end{remark}

According to the Weyl asymptotic formula
\begin{equation}\label{3.9}
N(\la)\ =\
C_n\,|\Om|\,\la^{n/2}+o(\la^{n/2})\,,
\qquad\la\to+\infty\,,
\end{equation}
(in the general case (\ref{3.9}) was proved in \cite{BS}).
The coefficient in the right hand side of (\ref{3.7})
contains an extra factor
$(1+2/n)^{n/2}$. G. P\'olya conjectured \cite{P}
that (\ref{3.7}) holds without this factor.
However, this remains an open problem.

Given a positive $\ve$, denote
\begin{equation*}
\Om_\ve^{\mathrm b}\ :=\ \{x\in\Om\,:\,d(x)\le\ve\}\,,
\qquad
\Om_\ve^{\mathrm i}\ :=\ \{x\in\Om\,:\,d(x)>\ve\}\,.
\end{equation*}
If
\begin{equation}\label{3.10}
|\Om_\ve^{\mathrm b}|\ \le\ \const\,\ve^r\,,\qquad r\in(0,1]\,,
\end{equation}
then, using the variational method \cite{CH},
one can prove that
\begin{equation}\label{3.11}
|\,N(\la)-C_n\,|\Om|\,\la^{n/2}\,|\ \le\
\begin{cases}
\const\,\la^{(n-1)/2}\ln\la\,,&r=1\,,\\
\const\,\la^{(n-r)/2}\,,&r<1\,.
\end{cases}
\end{equation}
It is well known that in the smooth case
$$
N(\la)-C_n\,|\Om|\,\la^{n/2}\ =\ O(\la^{(n-1)/2})
$$
(see, for example, \cite{I1} or \cite{SV}), but
it is not clear whether this estimate remains valid
for an arbitrary domain satisfying (\ref{3.10})
with $r=1$.

There is a number of papers devoted to estimates
of the remainder term in the Weyl formula.
In \cite{BL} the authors, applying the variational
technique, obtain explicit estimates for the constants
in (\ref{3.11}).
In order to prove the estimate of $N(\la)$ from
above, they imposed an additional condition on
the outer neighbourhood of the boundary $\dpa\Om$,
but this condition can probably be removed \cite{Ne}.
In \cite{Kr} the author estimated the remainder term
with the use of a different technique (similar to
that in \cite{LY}); his results seem to be less
precise than those obtained in \cite{BL}.

Let
$$
N_\ve^{\mathrm b}(\la):=
\int_{\Om_\ve^{\mathrm b}}e(x,x;\la)\,\dr x\,,
\qquad
N_\ve^{\mathrm i}(\la):=\int_{\Om_\ve^{\mathrm i}}
e(x,x;\la)\,\dr x\,.
$$
Then $\,N(\la)=N_\ve^{\mathrm b}(\la)+N_\ve^{\mathrm i}(\la)\,$
for every $\ve>0$.

\begin{corollary}\label{C3.4}
For all $\la>0$ and $\ve>0$ we have
\begin{eqnarray}
&&N_\ve^{\mathrm i}(\la)\ \ge\
C_n\,|\Om_\ve^{\mathrm i}|\,\la^{n/2}
-C_{n,1}\,(\la^{1/2}+\ve^{-1}\nu_{m_n})^{n-1}
\int_{\Om_\ve^{\mathrm i}}\frac{\dr x}{d(x)}\,,
\label{3.12}\\
&&N_\ve^{\mathrm i}(\la)\ \le\
C_n\,|\Om_\ve^{\mathrm i}|\,\la^{n/2}
+C_{n,2}\,(\la^{1/2}+\ve^{-1}\nu_{m_n})^{n-1}
\int_{\Om_\ve^{\mathrm i}}\frac{\dr x}{d(x)}\,,
\label{3.13}\\
&&N_\ve^{\mathrm b}(\la)\ \le\
C_{n,3}\,|\Om_\ve^{\mathrm b}|\,\la^{n/2}
+C_{n,4}\,\la^{-1/2}
(\la^{1/2}+\ve^{-1}\nu_{m_n})^{n-1}
\int_{\Om_\ve^{\mathrm i}}\frac{\dr x}{(d(x))^2}\,,
\label{3.14}\end{eqnarray}
where
$$
\begin{array}{lll}
&C_{n,1}=
n\,C_n\,2\pi^{-1}\nu_{m_n}^2\,,\qquad
&C_{n,2}=
n\,C_n\,(2\pi^{-1}\nu_{m_n}^2+\nu_{m_n})\,,\\
\\
&C_{n,3}=
(1+2/n)^{n/2}\,C_n\,,\qquad
&C_{n,4}=
(1+2/n)^{n/2}\,n^2\,C_n\,\nu_{m_n}^2\,.
\end{array}
$$
\end{corollary}

\begin{proof}
The inequalities (\ref{3.12}), (\ref{3.13}) are proved
by straightforward integration of (\ref{3.2}), (\ref{3.3}).
Theorem \ref{T3.2} and  (\ref{3.4}) imply that
\begin{multline}\label{3.15}
\int_0^\la N_\ve^{\mathrm b}(\mu)\,\dr\mu\ =\
\int_0^\la N(\la)\,\dr\mu\,-
\int_0^\la N_\ve^{\mathrm i}(\mu)\,\dr\mu\\
\le\ \frac{2}{n+2}\;C_n\,|\Om_\ve^{\mathrm b}|\,\la^{n/2+1}
+2n\,C_n\,\nu_{m_n}^2\,\la^{1/2}
\left(\la^{1/2}+\ve^{-1}\nu_{m_n}\right)^{n-1}
\int_{\Om_\ve^{\mathrm i}}\frac{\dr x}{(d(x))^2}\,.
\end{multline}
Now, applying the first inequality (\ref{3.8}) with
$\th=2/n$, we arrive at (\ref{3.14}).
\end{proof}

Adding up the inequalities (\ref{3.13}) and (\ref{3.14})
we obtain
\begin{multline}\label{3.16}
N(\la)\ \le\ C_n\,|\Om_\ve^{\mathrm i}|\,\la^{n/2}
+C_{n,3}\,|\Om_\ve^{\mathrm b}|\,\la^{n/2}\\
+\ (\la^{1/2}+\ve^{-1}\nu_{m_n})^{n-1}
\int_{\Om_\ve^{\mathrm i}}\frac{C_{n,2}\,d(x)
+C_{n,4}\,\la^{-1/2}}{(d(x))^2}\,\dr x\,,\qquad\forall\ve>0\,.
\end{multline}
Since
\begin{equation}\label{3.17}
\int_{\Om_\ve^{\mathrm i}}\frac{\dr x}{(d(x))^j}
\ =\ \int_\ve^\infty s^{-j}\,\dr(|\Om_s^{\mathrm b}|)
\ =\ j\int_\ve^\infty s^{-j-1}\,|\Om_s^{\mathrm b}|\,\dr s
-\,\ve^{-j}\,|\Om_\ve^{\mathrm b}|\,,
\end{equation}
(\ref{3.10}) and the inequalities (\ref{3.12}),
(\ref{3.16}) with $\ve=\la^{-1/2}$ imply (\ref{3.11}).

By (\ref{3.12}) and (\ref{3.16}) we have
\begin{eqnarray}
&&-\,C_n\,|\Om_\ve^{\mathrm b}|
\;-\;|\Om_\ve^{\mathrm i}|\;\frac{C_{n,1}}{\ve\la^{1/2}}\,
\left(1+\frac{\nu_{m_n}}{\ve\la^{1/2}}\right)^{n-1}\nonumber\\
&&\le\ \la^{-n/2}\,N(\la)\;-\;C_n\,|\Om|\label{3.18}\\
&&\le\ (C_{n,3}-C_n)\,|\Om_\ve^{\mathrm b}|
\;+\;|\Om_\ve^{\mathrm i}|\left(\frac{C_{n,2}}{\ve\la^{1/2}}\,
+\,\frac{C_{n,4}}{\ve^2\la}\right)
\left(1+\frac{\nu_{m_n}}{\ve\la^{1/2}}\right)^{n-1}\nonumber
\end{eqnarray}
for all $\ve>0$. If $\ve\to\infty$ then the second inequality
(\ref{3.18}) turns into (\ref{3.7}).
Since $|\Om_\ve^{\mathrm b}|\to0$ as $\ve\to0$, (\ref{3.18})
implies (\ref{3.9}). Moreover, taking $\ve=\la^{-\ka}$
with an arbitrary $\ka\in(0,\frac12)$, we
obtain the Weyl formula with a remainder estimate
\begin{equation*}
\la^{-n/2}N(\la)-C_n\,|\Om|\ =\
O(|\Om_{\la^{-\ka}}^{\mathrm b}|+\la^{\ka-1/2})\,,
\qquad\la\to+\infty\,.
\end{equation*}

\begin{remark}\label{R3.5}
If the condition (\ref{3.10}) is fulfilled then
integrating (\ref{3.4}) over $\Om_{\la^{-1/2}}^{\mathrm i}\,$,
applying (\ref{3.17}) and taking into account (\ref{3.6}),
we see that
\begin{equation}\label{3.19}
\la^{-1}\int_0^\la N(\mu)\,\dr\mu
\ =\
\frac{2}{n+2}\;C_n\,|\Om|\,\la^{n/2}+O(\la^{(n-r)/2})
\end{equation}
for all $r\in(0,1]$.
\end{remark}



\end{document}